\documentclass[reqno,report]{amsart}
\usepackage{latexsym}
\usepackage{amssymb}
\usepackage{amsmath}
\usepackage{amsfonts}
\usepackage{amsthm}
\usepackage{amscd}
\usepackage{texdraw}
\usepackage{color}
\usepackage[centertags]{amsmath}

\newlength{\defbaselineskip}
\newcommand{\setlinespacing}[1]%
           {\setlength{\baselineskip}{#1 \defbaselineskip}}

\theoremstyle{plain}
\newtheorem{thm}{Theorem}
\newtheorem{lem}{Lemma}
\newtheorem{prop}{Proposition}
\newtheorem{cor}{Corollary}

{\allowdisplaybreaks
\begin{document}
\title[Inversion of weighted Dirichlet series]{On inversion of absolutely convergent weighted Dirichlet series in two variables}

\author[P. A. Dabhi]{Prakash A. Dabhi}
\address{Department of Mathematics, Institute of Infrastructure Technology Research and Management(IITRAM), Maninagar (East), Ahmedabad - 380026, Gujarat, India}

\email{lightatinfinite@gmail.com, prakashdabhi@iitram.ac.in}
\subjclass{Primary 30B50; Secondary 30B50, 46H99, 46J05}
\thanks{The author is thankful to SERB, India, for the MATRICS grant no. MTR/2019/000162. The author is grateful to Professor Subhash J. Bhatt for suggesting this topic.}
\keywords{Dirichlet series, Fourier series, Wiener, L\'evy and \.Zelazko theorems, Beurling algebra, commutative Banach algebra, $p$-Banach algebra.}

\date{}

\dedicatory{}

\commby{}

\begin{abstract}
Let $0<p\leq 1$, and let $\omega:\mathbb N^2 \to [1,\infty)$ be an almost monotone weight. Let $\mathbb H$ be the closed right half plane in the complex plane. Let $\widetilde a$ be a complex valued function on $\mathbb H^2$ such that $\widetilde a(s_1,s_2)=\sum_{(m,n)\in \mathbb N^2}a(m,n)m^{-s_1}n^{-s_2}$ for all $(s_1,s_2)\in \mathbb H^2$ with $\sum_{(m,n)\in \mathbb N^2} |a(m,n)|^p\omega(m,n)<\infty$. If $\widetilde a$ is bounded away from zero on $\mathbb H^2$, then there is an almost monotone weight $\nu$ on $\mathbb N^2$ such that $1\leq \nu\leq \omega$, $\nu$ is constant if and only if $\omega$ is constant, $\nu$ is admissible if and only if $\omega$ is admissible, the reciprocal $\frac{1}{\widetilde a}$ has the Dirichlet representation $\frac{1}{\widetilde a}(s_1,s_2)=\sum_{(m,n)\in \mathbb N^2}b(m,n)m^{-s_1}n^{-s_2}$ for all $(s_1,s_2)\in \mathbb H^2$ and $\sum_{(m,n)\in \mathbb N^2}|b(m,n)|^p\nu(m,n)<\infty$. If $\varphi$ is holomorphic on a neighbourhood of the closure of range of $\widetilde a$, then there is an almost monotone weight $\xi$ on $\mathbb N^2$ such that $1\leq \xi\leq \omega$, $\xi$ is constant if and only if $\omega$ is constant, $\xi$ is admissible if and only if $\omega$ is admissible, $\varphi \circ \widetilde a$ has the Dirichlet series representation $(\varphi\circ \widetilde a)(s_1,s_2)=\sum_{(m,n)\in \mathbb N^2} c(m,n)m^{-s_1}n^{-s_2}\;((s_1,s_2)\in \mathbb H^2)$ and $\sum_{(m,n)\in \mathbb N^2}|c(m,n)|^p\xi(m,n)<\infty$. Let $\omega$ be an admissible weight on $\mathbb N^2$, and let $\widetilde a$ have $p$-th power $\omega$- absolutely convergent Dirichlet series. Then it is shown that the reciprocal of $\widetilde a$ has $p$-th power $\omega$- absolutely convergent Dirichlet series if and only if $\widetilde a$ is bounded away from zero. Though proofs are given for a function of two variables, it can naturally be extended for a function of several variables.
\end{abstract}

\maketitle

\section{Introduction}
If $f$ is a nowhere vanishing continuous function on the unit circle $\mathbb T=\{z\in \mathbb C:|z|=1 \}$ having absolutely convergent Fourier series, then the celebrated Wiener's theorem \cite{W} asserts that $\frac{1}{f}$ also has an absolutely convergent Fourier series. If $\varphi$ is holomorphic in a neighbourhood of the range of a continuous function $f$ having absolutely convergent Fourier series, then L\'evy's theorem \cite{L} implies that $\varphi\circ f$ has absolutely convergent Fourier series. The $p$-th power analogues of these theorems for $0<p\leq 1$ are obtained by \.Zelazko in \cite{Z}. A map $\omega:\mathbb Z \to [1,\infty)$ is a \emph{weight} if $\omega(m+n)\leq \omega(m)\omega(n)$ for all $m,n \in \mathbb Z$. Let $f$ be a continuous function on $\mathbb T$. Then the Fourier series of $f$ is \emph{$\omega$- absolutely convergent} if $\sum_{n\in \mathbb Z}|\widehat f(n)|\omega(n)<\infty$, where $\widehat f(n)=\frac{1}{2\pi}\int_0^{2\pi}f(e^{i\theta})e^{-in\theta}d\theta\;(n \in \mathbb Z)$ are the \emph{Fourier coefficients} of $f$. If $\omega$ is a Beurling-Domar weight, i.e., $\sum_{n\in \mathbb Z}\frac{\log \omega(n)}{1+n^2}<\infty$, $f$ has $\omega$- absolutely convergent Fourier series and if $f$ is nowhere zero on the circle $\mathbb T$, then Domar's Theorem \cite{D} implies that $\frac{1}{f}$ has $\omega$- absolutely convergent Fourier series. Weighted analogues and $p$-th power weighted analogues, $0<p\leq 1$, of theorems of Wiener, L\'evy and \.Zelazko are obtained in \cite{BDD,BhDe}. It is natural question to ask the same for Dirichlet series.

Let $\mathbb H=\{z\in \mathbb C:\text{Re }z\geq 0\}$ be the closed right half plane. A function $\widetilde a:\mathbb H\to \mathbb C$ has an \emph{absolutely convergent Dirichlet series} if $\widetilde a$ has the representation $$\widetilde a(s)=\sum_{n=1}^\infty a(n)n^{-s}\qquad(s \in \mathbb H),$$ and $\sum_{n=1}^\infty |a(n)|<\infty$. The above series is the \emph{Dirichlet series} of $\widetilde a$. Let $\widetilde a$ be a complex valued function on $\mathbb H$ having absolutely convergent Dirichlet series. It is shown in \cite{HW} that $\frac{1}{\widetilde a}$ has absolutely convergent Dirichlet series if and only if $\widetilde a$ is bounded away from zero on $\mathbb H$. A similar result for general Dirichlet series is obtained in \cite{Ed}. An elementary proof of the same result is given in \cite{GN}. Let $\omega$ be a weight on the multiplicative semigroup $\mathbb N$, i.e., $\omega\geq 1$ and $\omega(mn)\leq \omega(m)\omega(n)$ for all $m,n \in \mathbb N$. Then $\omega$ is an \emph{admissible weight} \cite{GL} if $\lim_{n\to \infty}\omega(k^n)^{\frac{1}{n}}=1$ for all $k \in \mathbb N$. Let $\omega$ be an admissible weight on $\mathbb N$, and let $\widetilde a$ have $\omega$- absolutely convergent Dirichlet series. It is recently shown in \cite{GL} that $\frac{1}{\widetilde a}$ has $\omega$- absolutely convergent Dirichlet series if and only if $\widetilde a$ is bounded away from zero.

We shall consider the semigroup $\mathbb N^2$ equipped with pointwise multiplication, i.e., $(m_1,n_1)(m_2,n_2)=(m_1m_2,n_1n_2)$ for all $(m_1,n_1), (m_2,n_2)\in \mathbb N^2$. Let $0<p\leq 1$, and let $\omega$ be a weight on $\mathbb N^2$, i.e., $$\omega((m_1,n_1)(m_2,n_2))\leq \omega(m_1,n_1)\omega(m_2,n_2)\quad ((m_1,n_1), (m_2,n_2)\in \mathbb N^2).$$ A function $\widetilde a:\mathbb H^2 \to \mathbb C$ has \emph{$p$-th power $\omega$- absolutely convergent Dirichlet series} if $\widetilde a$ has the representation $$\widetilde a(s_1,s_2)=\sum_{(m,n)\in \mathbb N^2}a(m,n)m^{-s_1}n^{-s_2}\;(s_1,s_2\in \mathbb H)$$ and $$\sum_{(m,n)\in \mathbb N^2}|a(m,n)|^p\omega(m,n)<\infty.$$ If $p=1$, we say $\widetilde a$ has $\omega$- absolutely convergent Dirichlet series.

Let $\omega$ be a weight on $\mathbb N^2$, and let $\{p_1,p_2,\ldots\}$ be the sequence of all primes arranged in increasing order. For $i \in \mathbb N$, define
\begin{eqnarray}\label{eq:1}
\rho_i=\lim_{n\to \infty}\omega(p_i^n,1)^{\frac{1}{n}}=\inf\{\omega(p_i^n,1)^{\frac{1}{n}}:n \in \mathbb N\}
\end{eqnarray}
and
\begin{eqnarray}\label{eq:2}
\mu_i=\lim_{n\to \infty}\omega(1,p_i^n)^{\frac{1}{n}}=\inf\{\omega(1,p_i^n)^{\frac{1}{n}}:n \in \mathbb N\}.
\end{eqnarray}
Then $\rho_i$ and $\mu_i$ are greater than or equal to $1$ for all $i\in \mathbb N$.

A weight $\omega$ on $\mathbb N^2$ is an \emph{admissible weight} if $\lim_{n\to \infty}\omega(l^n,k^n)^{\frac{1}{n}}=1$ for all $(l,k) \in \mathbb N^2$.

We call a weight $\omega$ on $\mathbb N^2$ to be an \emph{almost monotone} if either $\omega$ is admissible or there is $K>0$ such that $\omega(m_1,n_1)\leq K\omega(m_2,n_2)$ whenever $(m_1,n_1), (m_2,n_2) \in \mathbb N^2$ and either $m_1|m_2$ or $n_1|n_2$.

We prove the following Theorems.
\begin{thm}\label{thmA}
Let $\widetilde a:\mathbb H^2 \to \mathbb C$ have absolutely convergent Dirichlet series. Then the following statements are equivalent.
\begin{enumerate}
\item[(i)] $\frac{1}{\widetilde a}$ has absolutely convergent Dirichlet series.
\item[(ii)] $(a(m,n))$ has inverse in $\ell^1(\mathbb N^2)$.
\item[(iii)] $\widetilde a$ is bounded away from zero.
\end{enumerate}
\end{thm}
We use Theorem \ref{thmA} to prove the following analogue of theorems of Wiener and L\'evy in two variables.
\begin{thm}\label{thm:1}
Let $\omega$ be an almost monotone weight on $\mathbb N^2$, and let $\widetilde a:\mathbb H^2 \to \mathbb C$ have $\omega$- absolutely convergent Dirichlet series. Then the following statements hold.\\
(I) If $\widetilde a$ is bounded away from zero, then there is an almost monotone weight $\nu$ on $\mathbb N^2$ such that
\begin{enumerate}
\item $\nu$ is constant if and only if $\omega$ is constant;
\item $\nu$ is admissible if and only if $\omega$ is admissible;
\item $1\leq \nu \leq \omega$;
\item $\frac{1}{\widetilde a}$ has $\nu$- absolutely convergent Dirichlet series.
\end{enumerate}
(II) If $\varphi$ is holomorphic on the closure of the range of $\widetilde a$, then there is an almost monotone weight $\xi$ on $\mathbb N^2$ such that
\begin{enumerate}
\item $\xi$ is constant if and only if $\omega$ is constant;
\item $\xi$ is admissible if and only if $\omega$ is admissible;
\item $1\leq \xi \leq \omega$;
\item $\varphi \circ \widetilde a$ has $\xi$- absolutely convergent Dirichlet series.
\end{enumerate}
\end{thm}
\section{Proof of Theorem \ref{thmA}}
We first list the symbols we are going to use: $\mathbb D=\{z \in \mathbb C:|z|\leq 1\}$, $\mathbb D^\infty$ is the countable infinite product of $\mathbb D$ with itself, for $a \in \mathbb C$ and $r>0$, $B(a,r)=\{z\in \mathbb C:|z-a|<r\}$ and $\overline{B(a,r)}=\{z\in \mathbb C:|z-a|\leq r\}$, $\mathbb N_0=\mathbb N\cup\{0\}$, $\{p_1,p_2,\ldots\}$ is the set of all primes arranged in increasing order.

Let $(\mathcal B,\|\cdot\|)$ be a commutative Banach algebra. A nonzero linear functional $\varphi$ on $\mathcal B$ satisfying $\varphi(ab)=\varphi(a)\varphi(b)\;(a,b \in \mathcal B)$ is a \emph{complex homomorphism} on $\mathcal B$. Let $\Delta(\mathcal B)$ be the collection of all complex homomorphisms on $\mathcal B$. For $a \in \mathcal B$, let $\widehat a:\Delta(\mathcal B)\to \mathbb C$, $\widehat a(\varphi)=\varphi(a)\;(\varphi\in \Delta(\mathcal B))$, be the \emph{Gel'fand transform} of $a$. The weakest topology on $\Delta(\mathcal B)$ making each $\widehat a$ continuous is the \emph{Gel'fand topology} on $\Delta(\mathcal B)$. The set $\Delta(\mathcal B)$ with the Gel'fand topology is the \emph{Gel'fand space} of $\mathcal B$. If $a \in \mathcal B$, then $r(a)=\sup\{|\varphi(a)|:\varphi \in \Delta(\mathcal A)\}$ is the \emph{spectral radius} of $a$. We also notice that if $\mathcal B$ is a commutative unital Banach algebra and $a \in \mathcal B$, then $a$ is invertible in $\mathcal B$ if and only if $\widehat a(\varphi)=\varphi(a)\neq 0$ for all $\varphi \in \Delta(\mathcal B)$.

Let $$\ell^1(\mathbb N^2)=\{a:\mathbb N^2 \to \mathbb C:\|a\|=\sum_{(m,n)\in \mathbb N^2}|a(m,n)|<\infty\}.$$ Then $\ell^1(\mathbb N^2)$ is a commutative Banach algebra with the above norm and the convolution multiplication
\begin{eqnarray}\label{conv}
(a\star b)(m,n)=\sum_{u_1u_2=m\atop v_1v_2=n}a(u_1,v_1)b(u_2,v_2)\qquad(a, b \in \ell^1(\mathbb N^2), (m,n) \in \mathbb N^2).
\end{eqnarray}
We may write an element $a$ of $\ell^1(\mathbb N^2)$ as $a=\sum_{(m,n) \in \mathbb N^2}a(m,n)\delta_{(m,n)}$, where $\delta_{(m,n)}:\mathbb N^2 \to \mathbb C$ is defined by $\delta_{(m,n)}(m,n)=1$ and $\delta_{(m,n)}(l,k)=0$ if $(l,k)\neq (m,n)$. A \emph{semicharacter} on $\mathbb N^2$ is a nonzero map $\chi:\mathbb N^2\to \mathbb C$ satisfying $\chi((m,n)(l,k))=\chi(m,n)\chi(l,k)\;((m,n),(l,k) \in \mathbb N^2)$. Let $\widehat{\mathbb N^2}$ be the set of all bounded semicharacters on $\mathbb N^2$. Observe that if $\chi$ is a semicharacter on $\mathbb N^2$, then $\chi(1,1)=1$. We may identify $\widehat{\mathbb N^2}$ with the Gel'fand space $\Delta(\ell^1(\mathbb N^2))$ of $\ell^1(\mathbb N^2)$ via the map $\chi\mapsto \varphi_{\chi}$, where $$\varphi_\chi(a)=\sum_{(m,n)\in \mathbb N^2}a(m,n)\chi(m,n)\qquad(a \in \ell^1(\mathbb N^2,\omega)).$$

For a bounded semicharacter $\chi$ on $\mathbb N^2$, let $\chi_1(m)=\chi(m,1)$ and $\chi_2(n)=\chi(1,n)$ for all $m,n \in \mathbb N$. Then both $\chi_1$ and $\chi_2$ are bounded semicharacters on $\mathbb N$. Conversely, if $\chi_1$ and $\chi_2$ are bounded semicharacters on $\mathbb N$, then $\chi(m,n)=\chi_1(m)\chi_2(n)\;((m,n)\in \mathbb N^2)$ is a bounded semicharacter on $\mathbb N^2$. Since a semicharacter on $\mathbb N$ can be uniquely determined by its values at prime numbers, a semicharacter on $\mathbb N^2$ is uniquely determined by its values at the points of the form $(p,1)$ and $(1,p)$, where $p$ is prime.

Let $\mathbf A$ be the collection of all $\widetilde a:\mathbb H^2\to \mathbb C$ such that $$\widetilde a(s_1,s_2)=\sum_{(m,n)\in \mathbb N^2} a(m,n)m^{-s_1}n^{-s_2}\;(s_1,s_2\in \mathbb H)$$ and $$\|\widetilde a\|=\sum_{(m,n)\in \mathbb N^2} |a(m,n)|<\infty.$$
Then $\mathbf A$ is a unital commutative Banach algebra with the norm defined above and the usual multiplication of functions. Clearly, the  Banach algebras $\mathbf A$ and $\ell^1(\mathbb N^2)$ are isometrically isomorphic via the map $\widetilde a\mapsto (a(m,n))$.

For the following discussion, we follow the notations from \cite{HW}. Let $S=\{q_1,q_2,\ldots\}$ be a sequence of reals all greater than $1$ such that $\{\log q_1,\log q_2,\ldots\}$ is linearly independent over $\mathbb Q$, the field of rationals. For $\sigma\geq 0$, let $M_\sigma(S)$ be the collection of all semicharaters $\chi$ on $\mathbb N$ such that $|\chi(p_i)|=q_i^{-\sigma}$ for all $i\in \mathbb N$, and let $M_\infty$ consist the only semicharacter $\chi(1)=1$, $\chi(p_i)=0$ for all $i \in \mathbb N$. Let $L_\sigma(S)$ be the collection of all semicharacters on $\mathbb N$ such that $\chi(p_j)=q_j^{-\sigma-it}$ for all $j \in \mathbb N$, where $t$ is a real number that depends on $\chi$ but not on $j$. Clearly, both $M_\sigma(S)$ and $L_\sigma(S)$ are contained in the Gel'fand space $\Delta(\ell^1(\mathbb N))$ of $\ell^1(\mathbb N)$. We shall use the following two lemmas which are proved in \cite{HW}.

\begin{lem}\cite[Lemma 2]{HW}\label{HW2}
Let $q_1,q_2,\ldots,q_N$ be real numbers all greater than $1$. For $\sigma\geq 0$, let $$D_\sigma=\{(z_1,\ldots,z_N,w_1,\ldots,w_N)\in \mathbb C^{2N}:|z_i|\leq q_i^{-\sigma}, |w_i|\leq q_i^{-\sigma}, i=1,\ldots,N\}.$$ Let $\varphi$ be a complex polynomial in $2N$ variables, and let $\epsilon>0$. If  $|\varphi(z_1,\ldots,z_N,w_1,\ldots,w_N)|<\epsilon$ for some $(z_1,\ldots,z_N,w_1,\ldots,w_N) \in D_0$, then there exists $\sigma\geq 0$ and a point $(z_1',\ldots,z_N',w_1',\ldots,w_N')$ such that $|z_i'|=|w_i'|=q_i^{-\sigma}$ for all $i=1,\ldots,N$ and $|\varphi(z_1',\ldots,z_N',w_1',\ldots,w_N')|<\epsilon$.
\end{lem}

\begin{lem}\cite[Lemma 3]{HW}\label{HW3}
Let $\sigma\geq 0$. Then $L_\sigma(S)$ is dense in $M_\sigma(S)$ in the Gel'fand topology.
\end{lem}

\begin{lem}
Let $S=\{q_1,q_2,\ldots\}$ be a sequence of real numbers all greater than $1$ such that $\{\log q_1,\log q_2,\ldots\}$ is linearly independent over $\mathbb Q$. Let $\chi_1,\chi_2 \in M_\sigma(S)$, $N \in \mathbb N$ and $\epsilon>0$. Then there exist real numbers $t_1$ and $t_2$ such that $$\left|\prod_{r=1}^N\chi_1(p_r)^{\alpha_r}\prod_{r=1}^N  \chi_2(p_r)^{\beta_r}-\left(\prod_{r=1}^Nq_r^{\alpha_r}\right)^{-\sigma-it_1}\left(\prod_{r=1}^Nq_r^{\beta_r}\right)^{-\sigma-it_2}\right|<\epsilon\sum_{r=1}^N(\alpha_r+\beta_r)$$ whenever $\alpha_i,\beta_i\in \mathbb N_0$.
\end{lem}
\begin{proof}
As $\chi_1, \chi_2 \in M_\sigma(S)$, by Lemma \ref{HW3}, we get real numbers $t_1$ and $t_2$ such that $|\chi_j(p_r)-q_r^{-\sigma-it_j}|<\epsilon$ for all $r=1,2,\ldots,N$ and $j=1,2$. Let $m=p_1^{\alpha_1}\cdots p_N^{\alpha}$ and $n=p_1^{\beta_1}\cdots p_N^{\beta_N}$ be in $\mathbb N$, where $\alpha_i,\beta_i \in \mathbb N_0$. Then
\begin{eqnarray*}
&&|\chi_1(m)\chi_2(n)-m^{-\sigma-it_1}n^{-\sigma-it_2}|\\
& = & \left|\prod_{r=1}^N\chi_1(p_r)^{\alpha_r}\prod_{r=1}^N\chi_2(p_r)^{\beta_r}-\prod_{r=1}^N(q_r^{-\sigma-it_1})^{\alpha_r}
\prod_{r=1}^N(q_r^{-\sigma-it_2})^{\beta_r}\right|\\
& \leq & \left|\prod_{r=1}^N\chi_1(p_r)^{\alpha_r}\prod_{r=1}^N\chi_2(p_r)^{\beta_r}-
\prod_{r=1}^N(q_r^{-\sigma-it_1})^{\alpha_r}\prod_{r=1}^N\chi_2(p_r)^{\beta_r}\right|\\
&&+
\left|\prod_{r=1}^N(q_r^{-\sigma-it_1})^{\alpha_r}\prod_{r=1}^N\chi_2(p_r)^{\beta_r}
-\prod_{r=1}^N(q_r^{-\sigma-it_1})^{\alpha_r}
\prod_{r=1}^N(q_r^{-\sigma-it_2})^{\beta_r}\right|\\
& \leq & \left|\prod_{r=1}^N\chi_1(p_r)^{\alpha_r}-
\prod_{r=1}^N(q_r^{-\sigma-it_1})^{\alpha_r}\right|+\left|\prod_{r=1}^N\chi_2(p_r)^{\beta_r}-
\prod_{r=1}^N(q_r^{-\sigma-it_2})^{\beta_r}\right|\\
& < & (\alpha_1+\cdots+\alpha_N)\epsilon+(\beta_1+\cdots+\beta_N)\epsilon\\
& = & \epsilon\sum_{r=1}^N(\alpha_r+\beta_r).
\end{eqnarray*}
\end{proof}

\noindent
\textbf{Proof of Theorem \ref{thmA}:}
Since $\mathbf A$ and $\ell^1(\mathbb N^2)$ are isometrically isomorphic as Banach algebras, the statements (i) and (ii) are equivalent. Assume that $\widetilde a$ has inverse in $\mathbf {A}$, i.e., $\frac{1}{\widetilde a}$ has absolutely convergent Dirichlet series. Then $\widehat{\widetilde a}(\chi)=\chi(\widetilde a)\neq 0$ for all $\chi \in \Delta(\mathbf{A})$. As $\widehat{\widetilde a}$ is continuous and $\Delta(\mathbf A)$ is compact, there is $k>0$ such that $|\widehat{\widetilde a}(\chi)|\geq k$ for all $\chi \in \Delta(\mathbf A)$. For $(s_1,s_2)\in \mathbb H^2$, let $$\chi_{s_1,s_2}(\widetilde b)=\sum_{(m,n)\in \mathbb N^2}b(m,n)m^{-s_1}n^{-s_2}\qquad(\widetilde b \in \mathbf{A}).$$ Then $\chi_{s_1,s_2} \in \Delta(\mathbf{A})$. It follows that $|\widetilde a(s_1,s_2)|=|\widehat{\widetilde a}(\chi_{s_1,s_2})|\geq k$ for all $(s_1,s_2)\in \mathbb H^2$. So, (i) $\Rightarrow$ (iii) is proved.

We now prove that (iii) $\Rightarrow$ (i). The proof is inspired by the proof of Theorem 1 in \cite{HW}. Suppose that $\widetilde a$ has no inverse in $\mathbf{A}$. Then there is a bounded semicharacter $\chi$ on $\mathbb N^2$ such that $\chi(\widetilde a)=0$. Since $\chi$ is a bounded semicharacter on $\mathbb N^2$, there exist bounded semicharacters $\chi_1$ and $\chi_2$ on $\mathbb N$ such that $\chi(m,n)=\chi_1(m)\chi_2(n)$ for all $(m,n)\in \mathbb N^2$. For $M \in \mathbb N$, let $\Box_M=\{1,2,\ldots,M\}\times \{1,2,\ldots,M\}$. Take any $\epsilon>0$. As $\sum_{(m,n)\in \mathbb N^2}|a(m,n)|<\infty$, there is $M \in \mathbb N$ such that $\sum_{(m,n)\in \mathbb N^2\setminus \Box_M}|a(m,n)|<\epsilon$. Since $\chi(\widetilde a)=0$, i.e., $\sum_{(m,n)\in \mathbb N^2}a(m,n)\chi(m,n)=0$ and $|\chi(m,n)|\leq 1$ for all $(m,n)$,
\begin{eqnarray*}
\left|\sum_{(m,n)\in \Box_M}a(m,n)\chi(m,n)\right| &=& \left|\sum_{(m,n)\in \mathbb N^2\setminus \Box_M}a(m,n)\chi(m,n)\right|\\
&\leq &\sum_{(m,n)\in \mathbb N^2\setminus \Box_M}|a(m,n)|<\epsilon.
\end{eqnarray*}
Observe that $\sum_{(m,n)\in \Box_M}a(m,n)\chi(m,n)=\sum_{(m,n)\in \Box_M}a(m,n)\chi_1(m)\chi_2(n)$ is a polynomial in $\chi_1(p_1),\ldots,\chi_1(p_N),\chi_2(p_1),\ldots,\chi_2(p_N)$, say $\varphi$, where $p_N$ is the largest prime less than or equal to $M$, i.e,
\begin{eqnarray*}
\varphi(\chi_1(p_1),\ldots,\chi_1(p_N),\chi_2(p_1),\ldots,\chi_2(p_N))&=&\sum_{(m,n)\in \Box_M}a(m,n)\chi(m,n)\\
&=&\sum_{(m,n)\in \Box_M}a(m,n)\chi_1(m)\chi_2(n).
\end{eqnarray*}
By Lemma \ref{HW2}, with $q_i=p_i$, $i=1,2,\ldots,N$, we get $\sigma\geq 0$ and an element $(z_1',\ldots,z_N',w_1',\ldots,w_N')\in \mathbb C^{2N}$ such that $|\varphi(z_1',\ldots,z_N',w_1',\ldots,w_N')|<\epsilon$, where $|z_i'|=|w_i'|=p_i^{-\sigma}$, $i=1,2,\ldots,N$. Define $\psi_1,\psi_2:\mathbb N \to \mathbb C$ by $\psi_1(p_i)=z_i'$, $\psi_2(p_i)=w_i'$, $i=1,2,\ldots,N$ and $\psi_1(p_i)=\psi_2(p_i)=p_i^{-\sigma}$ otherwise. Then $\psi_1$ and $\psi_2$ are bounded semicharacters on $\mathbb N$ and hence $\psi=\psi_1\psi_2$ is a bounded semicharacter on $\mathbb N^2$. As $\psi_1, \psi_2 \in M_\sigma(S)$, by Lemma \ref{HW3} with $S=\{p_1,p_2,\ldots\}$, we get real numbers $t_1$ and $t_2$ such that
\begin{eqnarray*}
\left|\prod_{r=1}^N\psi_1(p_r)^{\alpha_r}\prod_{r=1}^N  \psi_2(p_r)^{\beta_r}-\left(\prod_{r=1}^Np_r^{\alpha_r}\right)^{-\sigma-it_1}\left(\prod_{r=1}^Np_r^{\beta_r}\right)^{-\sigma-it_2}\right|
<\frac{\epsilon}{QJ}\sum_{r=1}^N(\alpha_r+\beta_r),
\end{eqnarray*}
whenever $\alpha_i,\beta_i \in \mathbb N_0$, where $Q$ is the degree of the polynomial $\varphi$ and $J=\sum_{(m,n)\in \Box_M}|a(m,n)|$. So, if $m=p_1^{\alpha_1}\cdots p_N^{\alpha_N}$ and $n=p_1^{\beta_1}\cdots p_N^{\beta_N}$ are less than or equal to $M$, where $\alpha_i,\beta_i \in \mathbb N_0$, then $|\psi(m,n)-m^{-\sigma-it_1}n^{-\sigma-it_2}|<\frac{\epsilon}{J}$ and hence
\begin{eqnarray*}
&&\left|\sum_{(m,n)\in\Box_M}a(m,n)\psi(m,n)-
\sum_{(m,n)\in\Box_M}a(m,n)m^{-\sigma-it_1}n^{-\sigma-it_2}\right|\\
& \leq & \sum_{(m,n)\in\Box_M}|a(m,n)||\psi(m,n)-m^{-\sigma-it_1}n^{-\sigma-it_2}|\\
&< & \sum_{(m,n)\in\Box_M}|a(m,n)|\frac{\epsilon}{J}=\epsilon.
\end{eqnarray*}
Now,
\begin{eqnarray*}
|\widetilde a(\sigma+it_1,\sigma+it_2)|&=&\left|\sum_{(m,n)\in\mathbb N^2}a(m,n)m^{-\sigma-it_1}n^{-\sigma-it_2}\right|\\
& \leq &\left|\sum_{(m,n)\in\Box_M}a(m,n)m^{-\sigma-it_1}n^{-\sigma-it_2}\right|+\sum_{(m,n)\in\mathbb N^2\setminus\Box_M}|a(m,n)|\\
& \leq & \left|\sum_{(m,n)\in\Box_M}a(m,n)m^{-\sigma-it_1}n^{-\sigma-it_2}-\sum_{(m,n)\in\Box_M}a(m,n)\psi(m,n)\right|\\
&& +\left|\sum_{(m,n)\in\Box_M}a(m,n)\psi(m,n)\right|+\epsilon\\
&<&3\epsilon.
\end{eqnarray*}
Thus, $\widetilde a$ is not bounded away from zero. $\hfill \Box$

\section{Proof of Theorem \ref{thm:1}}
Let $\omega:\mathbb N^2 \to [1,\infty)$  be a weight. Then $$\ell^1(\mathbb N^2,\omega)=\{a:\mathbb N^2 \to \mathbb C:\|a\|=\sum_{(m,n)\in \mathbb N^2}|a(m,n)|\omega(m,n)<\infty\}$$ is a commutative Banach algebra with the above norm and the convolution multiplication defined by (\ref{conv}).

Let $\mathbf A_\omega$ be the collection of all $\widetilde a:\mathbb H^2\to \mathbb C$ such that $$\widetilde a(s_1,s_2)=\sum_{(m,n)\in \mathbb N^2} a(m,n)m^{-s_1}n^{-s_2}\;(s_1,s_2\in \mathbb H)$$ and $$\|\widetilde a\|=\sum_{(m,n)\in \mathbb N^2} |a(m,n)|\omega(m,n)<\infty.$$
When $\omega=1$, we write $\mathbf A$ in place of $\mathbf A_\omega$. Then $\mathbf A_\omega$ is a unital commutative Banach algebra with the norm defined above and the usual multiplication of functions and is isometrically isomorphic to $\ell^1(\mathbb N^2,\omega)$ via the map $\widetilde a\mapsto (a(m,n))$.

A semicharacter $\chi$ on $\mathbb N^2$ is \emph{$\omega$- bounded} if $|\chi(m,n)|\leq \omega(m,n)\;((m,n)\in \mathbb N^2)$. Let $\widehat{\mathbb N^2_\omega}$ be the set of all $\omega$- bounded semicharacters on $\mathbb N^2$. In the weighted case, $\widehat{\mathbb N^2_\omega}$ can be idenfied with the Gel'fand space $\Delta(\ell^1(\mathbb N^2,\omega))$ of $\ell^1(\mathbb N^2,\omega)$ via the map $\chi\mapsto \varphi_{\chi}$, where $$\varphi_\chi(a)=\sum_{(m,n)\in \mathbb N^2}a(m,n)\chi(m,n)\qquad(a \in \ell^1(\mathbb N^2,\omega)).$$

\begin{lem}\label{lem:1}
Let $\omega$ be a weight on $\mathbb N^2$. Then the Gel'fand space of $\ell^1(\mathbb N^2,\omega)$ can be identified with the collection of all pairs $((s_i),(t_i))$ of complex sequences satisfying
\begin{eqnarray}\label{semch}
|s_{i_1}|^{\alpha_1}\cdots|s_{i_k}|^{\alpha_k}|t_{i_1}|^{\beta_1}\cdots|t_{i_k}|^{\beta_k}\leq \omega(p_{i_1}^{\alpha_1}\cdots p_{i_k}^{\alpha_k},p_{i_1}^{\beta_1}\cdots p_{i_k}^{\beta_k}),
\end{eqnarray}whenever $p_{i_1},p_{i_2}, \ldots,p_{i_k}$ are distinct primes, $\alpha_i,\beta_i \in \mathbb N_0$ and $k \in \mathbb N$.
\end{lem}
\begin{proof}
Since $\widehat{\mathbb N^2_\omega}$ can be identified with the Gel'fand space $\Delta(\ell^1(\mathbb N^2,\omega))$ of $\ell^1(\mathbb N^2,\omega)$, we shall identify $\widehat{\mathbb N^2_\omega}$ with the collection of all pair of complex sequences satisfying (\ref{semch}).

Let $\chi \in \widehat{\mathbb N^2_\omega}$. Let $s_i=\chi(p_i,1)$ and $t_i=\chi(1,p_i)$ for all $i\in \mathbb N$. Let $p_{i_1},\ldots,p_{i_k}$ be distinct primes, and let $\alpha_i,\beta_i \in \mathbb N_0$. Since $\chi$ is $\omega$- bounded, we get
\begin{eqnarray*}
&&|s_{i_1}|^{\alpha_1}\cdots|s_{i_k}|^{\alpha_k}|t_{i_1}|^{\beta_1}\cdots|t_{i_k}|^{\beta_k}\\
& = & |\chi(p_{i_1},1)|^{\alpha_1}\cdots|\chi(p_{i_k},1)|^{\alpha_k}|
\chi(1,p_{i_1})|^{\beta_1}\cdots|\chi(1,p_{i_1})|^{\beta_k}\\
& =& |\chi(p_{i_1}^{\alpha_1}\cdots p_{i_k}^{\alpha_k},p_{i_1}^{\beta_1}\cdots p_{i_k}^{\beta_k})|\\
& \leq & \omega(p_{i_1}^{\alpha_1}\cdots p_{i_k}^{\alpha_k},p_{i_1}^{\beta_1}\cdots p_{i_k}^{\beta_k}).
\end{eqnarray*}
Conversely, let $((s_i),(t_i))$ be a pair of complex sequences satisfying (\ref{semch}). Define $\chi:\mathbb N^2\to \mathbb C$ by $\chi(1,1)=1$, $\chi(p_i,1)=s_i$ and $\chi(1,p_i)=t_i$ for all $i \in \mathbb N$ and extend on $\mathbb N^2$ so that it is multiplicative. Then $\chi$ is a semicharacter on $\mathbb N^2$. Let $m=p_{i_1}^{\alpha_1}\cdots p_{i_j}^{\alpha_j}$ and $n=p_{i_1}^{\beta_1}\cdots p_{i_j}^{\beta_j}$ be in $\mathbb N$, where $p_{i_1},\ldots,p_{i_j}$ are distinct primes and $\alpha_i,\beta_i \in \mathbb N_0$. Then
\begin{eqnarray*}
|\chi(m,n)| & = & |\chi(p_{i_1}^{\alpha_1}\cdots p_{i_j}^{\alpha_j},p_{i_1}^{\beta_1}\cdots p_{i_j}^{\beta_j})\\
& = &|\chi(p_{i_1},1)|^{\alpha_1}\cdots |\chi(p_{i_j},1)|^{\alpha_j}|\chi(1,p_{i_1})|^{\beta_1}\cdots |\chi(1,p_{i_j})|^{\beta_j}\\
&=& |s_{i_1}|^{\alpha_1}\cdots |s_{i_j}|^{\alpha_j}|t_{i_1}|^{\beta_1}\cdots |t_{i_j}|^{\beta_j}\\
& \leq & \omega(p_{i_1}^{\alpha_1}\cdots p_{i_j}^{\alpha_j},p_{i_1}^{\beta_1}\cdots p_{i_j}^{\beta_j})\\
& = & \omega(m,n).
\end{eqnarray*}
Therefore $\chi$ is an $\omega$- bounded semicharacter on $\mathbb N^2$.
\end{proof}
Using the above lemma, an $\omega$- bounded semicharacter $\chi$ on $\mathbb N^2$ can be identified with $((s_i),(t_i))$ satisfying the inequality (\ref{semch}). We write this $\chi$ by $\chi_{((s_i),(t_i))}$. We also notice that if $\chi\in \widehat{\mathbb N^2_\omega}$, then $|\chi(p_i,1)|\leq \rho_i$ and $|\chi(1,p_i)|\leq \mu_i$ for all $i$.

\begin{cor}\label{cor:1}
Let $\omega$ be a weight on $\mathbb N^2$, and let $(r_i)$ and $(s_i)$ be sequences of positive real numbers such that $$r_{i_1}^{\alpha_1}r_{i_2}^{\alpha_2}\cdots r_{i_k}^{\alpha_k}s_{i_1}^{\beta_1}s_{i_2}^{\beta_2}\cdots s_{i_k}^{\beta_k}\leq \omega(p_{i_1}^{\alpha_1}p_{i_2}^{\alpha_2}\cdots p_{i_k}^{\alpha_k},p_{i_1}^{\beta_1}p_{i_2}^{\beta_2}\cdots p_{i_k}^{\beta_k})$$ whenever $p_{i_j}$ are distinct primes, $\alpha_i,\beta_i \in \mathbb N_0$ and $k \in \mathbb N$. Then $\prod_{i=1}^\infty \overline{B(0,r_i)}\times \prod_{i=1}^\infty \overline{B(0,s_i)}\subset \Delta(\ell^1(\mathbb N^2,\omega))\subset \prod_{i=1}^\infty \overline{B(0,\rho_i)}\times \prod_{i=1}^\infty \overline{B(0,\mu_i)}$.
\end{cor}

\begin{lem}\label{lem:2}
Let $\omega$ be a weight on $\mathbb N^2$. Then $\omega$ is an admissible weight if and only if $\rho_i=\mu_i=1$ for all $i \in \mathbb N$.
\end{lem}
\begin{proof}
Clearly, if $\omega$ is an admissible weight, then $\lim_{n\to \infty}\omega(p^n,1)^{\frac{1}{n}}=1$ and $\lim_{n\to \infty}\omega(1,p^n)^{\frac{1}{n}}=1$ for every prime $p$. So, $\rho_i=\mu_i=1$ for all $i \in \mathbb N$.

Conversely, assume that $\rho_i=\mu_i=1$ for all $i \in \mathbb N$. Let $l,k\in \mathbb N$. Then we get distinct primes $p_1,p_2,\ldots,p_m$ and $\alpha_1,\alpha_2,\ldots,\alpha_m,\beta_1,\beta_2,\ldots,\beta_m \in \mathbb N_0$ such that $l=p_1^{\alpha_1}\cdots p_m^{\alpha_m}$ and $k=p_1^{\beta_1}\cdots p_m^{\beta_m}$. Let $n \in \mathbb N$. Then
\begin{eqnarray*}
1 &\leq & \omega(l^n,k^n)^{\frac{1}{n}}=\omega(p_1^{\alpha_1 n}\cdots p_m^{\alpha_m n},p_1^{\beta_1 n}\cdots p_m^{\beta_m n})^{\frac{1}{n}}\\
&\leq & \omega(p_1^n,1)^{\frac{\alpha_1}{n}}\cdots \omega(p_m^n,1)^{\frac{\alpha_m}{n}}\omega(1,p_1^n)^{\frac{\beta_1}{n}}\cdots \omega(1,p_m^n)^{\frac{\beta_m}{n}}.
\end{eqnarray*}
It follows that $\lim_{n\to \infty}\omega(l^n,k^n)^{\frac{1}{n}}=1$.
\end{proof}

\begin{cor}\label{cor:3}
Let $\omega$ be a weight on $\mathbb N^2$. Then $\omega$ is an admissible weight if and only if $\Delta(\ell^1(\mathbb N^2,\omega))=(\mathbb D^\infty)^2$.
\end{cor}
\begin{proof}
Assume that $\omega$ is admissible. Since $\omega\geq 1$, $(\mathbb D^\infty)^2 \subset \Delta(\ell^1(\mathbb N^2,\omega))$. Let $\chi\in \Delta(\ell^1(\mathbb N^2,\omega))$, and let $p$ be any prime. Then $|\chi(p,1)|=|\chi(p^n,1)|^{\frac{1}{n}}\leq \omega(p^n,1)^{\frac{1}{n}}$ for all $n \in \mathbb N$. This implies that $|\chi(p,1)|\leq 1$. Similarly, $|\chi(1,p)|\leq 1$. Therefore $\Delta(\ell^1(\mathbb N^2,\omega))\subset (\mathbb D^\infty)^2$.

Conversely, assume that $(\mathbb D^\infty)^2= \Delta(\ell^1(\mathbb N^2,\omega))$. Let $p$ be a prime. By assumption $|\chi(p,1)|\leq 1$ and $|\chi(1,p)|\leq 1$ for all $\chi \in \widehat{\mathbb N}_\omega$. Now,
\begin{eqnarray*}
1 & \leq & \lim_{n\to \infty}\omega(p^n,1)^{\frac{1}{n}}\\
& = & r(\delta_{(p,1)})=\sup\{|\varphi_\chi(\delta_{(p,1)})|:\chi \in \widehat{\mathbb N^2_\omega}\}\\
& = &\sup\{|\chi(p,1)|:\chi \in \widehat{\mathbb N^2_\omega}\}\leq 1.
\end{eqnarray*}
Also, $\lim_{n\to \infty}\omega(1,p^n)^{\frac{1}{n}}=1$. By Lemma \ref{lem:2}, $\omega$ is an admissible weight.
\end{proof}

\begin{cor}\label{cor:2}
Let $\omega$ be an almost monotone weight on $\mathbb N^2$, and let $$X_j=(\mathbb D\times \mathbb D\times \cdots \mathbb D\times \overline{B(0,\rho_j)}\times \mathbb D\times \cdots)\times \mathbb D^\infty,$$ where $\overline{B(0,\rho_j)}$ is at the $j$-th place. Then $X_j \subset \Delta(\ell^1(\mathbb N^2,\omega))$ for all $j \in \mathbb N$.
\end{cor}
\begin{proof}
If $\omega$ is admissible, then $X_j=(\mathbb D^\infty)^2=\Delta(\ell^1(\mathbb N^2,\omega))$. Now, assume that $\omega$ is not admissible. As $\omega$ is an almost monotone weight, there is $K>0$ such that $\omega(m_1,n_1)\leq K\omega(m_2,n_2)$ whenever $(m_1,n_1), (m_2,n_2) \in \mathbb N^2$ and $m_1|m_2$ or $n_1|n_2$. Let $((s_i),(t_i))\in X_j$, i.e., $\chi_{((s_i),(t_i))}(p_k,1)=s_k$ and $\chi_{((s_i),(t_i))}(1,p_k)=t_k$ for all $k\in \mathbb N$ and $\chi_{((s_i),(t_i))}(1,1)=1$. Let $m=p_{i_i}^{\alpha_1}\cdots p_{i_k}^{\alpha_k}$ and $n=p_{l_1}^{\beta_1}\cdots p_{l_u}^{\beta_u}$, where $p_{i_1},\ldots,p_{i_j}$ are distinct primes, $p_{l_1},\ldots,p_{l_u}$ are distinct primes and $\alpha_i,\beta_i \in \mathbb N$. Assume that $p_j$ is not equal to any of $p_{i_1},\ldots,p_{i_j}$. Then
\begin{eqnarray*}
|\chi_{((s_i),(t_i))}(m,n)| & = & |\chi_{((s_i),(t_i))}(p_{i_i}^{\alpha_1}\cdots p_{i_k}^{\alpha_k},p_{l_1}^{\beta_1}\cdots p_{l_u}^{\beta_u})|\\
& = & |s_{i_1}|^{\alpha_1}\cdots |s_{i_j}|^{\alpha_j}|t_{l_1}|^{\beta_1}\cdots |t_{l_u}|^{\beta_u}\\
& \leq & 1 \leq \omega(m,n).
\end{eqnarray*}
Now, assume that $p_j$ is one of $p_{i_1},\ldots,p_{i_j}$. We may assume that $p_j=p_{i_1}$. Then
\begin{eqnarray*}
|\chi_{((s_i),(t_i))}(m,n)| & = & |\chi_{((s_i),(t_i))}(p_{i_i}^{\alpha_1}\cdots p_{i_k}^{\alpha_k},p_{l_1}^{\beta_1}\cdots p_{l_u}^{\beta_u})|\\
& = & \rho_j^{\alpha_1}|s_{i_2}|^{\alpha_2}\cdots |s_{i_j}|^{\alpha_j}|t_{l_1}|^{\beta_1}\cdots |t_{l_u}|^{\beta_u}\\
& \leq & \rho_j^{\alpha_1}\\
& \leq & \omega(p_j^{\alpha_1},1)\\
&\leq & K\omega(p_{i_i}^{\alpha_1}\cdots p_{i_k}^{\alpha_k},p_{l_1}^{\beta_1}\cdots p_{l_u}^{\beta_u})\\
&=& K\omega(m,n).
\end{eqnarray*}
If $k \in \mathbb N$, then by using the same arguments above, we have $|\chi_{((s_i),(t_i))}(m,n)|^k=|\chi_{((s_i),(t_i))}(m^k,n^k)|\leq K\omega(m^k,n^k)\leq K\omega(m,n)^k$. This gives $|\chi_{((s_i),(t_i))}(m,n)|\leq \omega(m,n)$. Thus $X_j \subset \Delta(\ell^1(\mathbb N,\omega))$.
\end{proof}

\begin{lem}\label{lem:3}
Let $(R_i)$ and $(S_i)$ be sequences of real numbers, all greater than or equal to $1$. Define $\nu$ on $\mathbb N^2$ as $\nu(1,1)=1$, $\nu(p_i,1)=R_i$ and $\nu(1,p_i)=S_i$ for $i \in \mathbb N$ and extend it on $\mathbb N^2$ so that it is multiplicative. Then $\nu$ is an almost monotone weight on $\mathbb N^2$ and $\Delta(\ell^1(\mathbb N^2,\nu))=\prod_{i=1}^\infty \overline{B(0,R_i)} \times \prod_{i=1}^\infty \overline{B(0,S_i)}$.
\end{lem}
\begin{proof}
Clearly, $\nu$ is an almost monotone weight. Let $P=\prod_{i=1}^\infty\overline{B(0,R_i)}\times \prod_{i=1}^\infty \overline{B(0,S_i)}$. Let $((\delta_i),(\eta_i)) \in P$. Then $|\delta_i|\leq R_i$ and $|\eta_i|\leq S_i$ for all $i\in \mathbb N$. Define $\chi:\mathbb N^2 \to \mathbb C$ by $\chi(1,1)=1$, $\chi(p_i,1)=\delta_i$ and $\chi(1,p_i)=\eta_i$ for all $i$ and extend it on $\mathbb N^2$ so that it is multiplicative. Then $\chi$ is a semicharacter on $\mathbb N^2$. Let $m=p_{i_1}^{\alpha_1}p_{i_2}^{\alpha_2}\cdots p_{i_k}^{\alpha_k}$ and $n=p_{i_1}^{\beta_1}p_{i_2}^{\beta_2}\cdots p_{i_k}^{\beta_k}$ be elements of $\mathbb N$, where $p_{i_j}$ are distinct primes and $\alpha_i,\beta_i \in \mathbb N_0$. Then
\begin{eqnarray*}
|\chi(m,n)|& = &|\chi(p_{i_1}^{\alpha_1}p_{i_2}^{\alpha_2}\cdots p_{i_k}^{\alpha_k},p_{i_1}^{\beta_1}p_{i_2}^{\beta_2}\cdots p_{i_k}^{\beta_k})|\\
& = & |\chi(p_{i_1},1)|^{\alpha_1}|\cdots |\chi(p_{i_k},1)|^{\alpha_k}|\chi(1,p_{i_1})|^{\beta_1}|\cdots |\chi(1,p_{i_k})|^{\beta_k}\\
& = & |\delta_{i_1}|^{\alpha_1}|\delta_{i_2}|^{\alpha_2}\cdots |\delta_{i_k}|^{\alpha_k}|\eta_{i_1}|^{\beta_1}|\eta_{i_2}|^{\beta_2}\cdots |\eta_{i_k}|^{\beta_k}\\
& \leq & R_{i_1}^{\alpha_1}R_{i_2}^{\alpha_2}\cdots R_{i_k}^{\alpha_k}S_{i_1}^{\beta_1}S_{i_2}^{\beta_2}\cdots S_{i_k}^{\beta_k}\\
& = & \nu(p_{i_1},1)^{\alpha_1}\nu(p_{i_2},1)^{\alpha_2}\cdots \nu(p_{\i_k},1)^{\alpha_k}\nu(1,p_{i_1})^{\beta_1}\nu(1,p_{i_2})^{\beta_2}\cdots \nu(1,p_{\i_k})^{\beta_k}\\
& = & \nu(p_{i_1}^{\alpha_1}p_{i_2}^{\alpha_2}\cdots p_{i_k}^{\alpha_k},p_{i_1}^{\beta_1}p_{i_2}^{\beta_2}\cdots p_{i_k}^{\beta_k})\\
&=&\nu(m,n).
\end{eqnarray*}
Therefore $\chi$ is a $\nu$- bounded semicharacter on $\mathbb N^2$. This proves that $P\subset \Delta(\ell^1(\mathbb N^2,\nu))$.

Let $\chi \in \Delta(\ell^1(\mathbb N^2,\nu))$. Then $\chi(1,1)=1$. For $i \in \mathbb N$, let $\delta_i=\chi(p_i,1)$ and $\eta_i=\chi(1,p_i)$. Then $|\delta_i|=|\chi(p_i,1)|\leq \nu(p_i,1)=R_i$ and $|\eta_i|=|\chi(1,p_i)|\leq \nu(1,p_i)=S_i$ for all $i\in \mathbb N$. So, $\chi \in P$. Therefore $\Delta(\ell^1(\mathbb N^2,\nu))\subset P$.
\end{proof}

\begin{prop}\label{prop:1}
Let $\omega$ be a weight on $\mathbb N^2$. Then the product topology on $(\mathbb D^\infty)^2$ is the relative Gel'fand topology on $(\mathbb D^\infty)^2$ as a subspace of $\Delta(\ell^1(\mathbb N^2,\omega))$.
\end{prop}
\begin{proof}
As $\omega\geq 1$, $(\mathbb D^\infty)^2\subset \Delta(\ell^1(\mathbb N^2,\omega))$. We first show that the product topology on $(\mathbb D^\infty)^2$ is weaker than the relative Gel'fand topology on it. Take a subbasis element of the form $V=(B(a,\epsilon)\times \mathbb D\times \cdots)\times \mathbb D^\infty$, where $a \in \mathbb D$ and $\epsilon>0$. Here $B(a,\epsilon)$ is the open ball with centre $a$ and radius $\epsilon$ in $\mathbb D$. Let $\chi:\mathbb N^2\to \mathbb C$ be defined by $\chi(1,1)=1$, $\chi(p_1,1)=a$, $\chi(p_i,1)=1$ for all $i\geq 2$ and $\chi(1,p_i)=1$ for all $i$. Let $w=((a,1,1,\ldots),(1,1,1,\ldots))$. Then $\chi=\chi_w$ is an $\omega$- bounded semicharacter on $\mathbb N^2$, i.e., $\chi_w \in \Delta(\ell^1(\mathbb N^2,\omega))$. The set $$U(\chi_w,\delta_{(p_1,1)},\epsilon)=\{\chi_{((z_i),(\eta_i))}\in \Delta(\ell^1(\mathbb N^2,\omega)):|\chi_{((z_i),(\eta_i))}(\delta_{(p_1,1)})-\chi_w(\delta_{(p_1,1)})|<\epsilon\}$$ is open in the Gel'fand topology on $\Delta(\ell^1(\mathbb N^2,\omega))$. But $U(\chi_w,\delta_{(p_1,1)},\epsilon)=\{\chi_{((z_i),(\eta_i))}\in \Delta(\ell^1(\mathbb N^2,\omega)):|z_1-a|<\epsilon\}$. Therefore $V=U(\chi_w,\delta_{(p_1,1)},\epsilon)\cap(\mathbb D^\infty)^2$ and hence the product topology on $(\mathbb D^\infty)^2$ is weaker than the relative Gel'fand topology on it.

We now show that $(\mathbb D^\infty)^2$ is closed in $\Delta(\ell^1(\mathbb N^2,\omega))$. Let $\chi_{((z_i),(\eta_i))}\in \Delta(\ell^1(\mathbb N^2,\omega))\setminus (\mathbb D^\infty)^2$. Then $|z_i|>1$ or $|\eta_i|>1$ for some $i$. We may assume that $|z_1|>1$. Then there is $\epsilon>0$ such that $B(z_1,\epsilon)\cap \mathbb D=\emptyset$. Note that $U(\chi_{((z_i),(\eta_i))},\delta_{(p_1,1)},\epsilon)$ is an open subset of $\Delta(\ell^1(\mathbb N^2,\omega))$ and it contains $\chi_{((z_i),(\eta_i))}$. Suppose that $\chi_{((\alpha_i),(\beta_i))}\in U(\chi_{((z_i),(\eta_i))},\delta_{(p_1,1)},\epsilon) \cap (\mathbb D^\infty)^2$. Then $|z_1-\alpha_1|<\epsilon$, so, $\alpha_1 \notin \mathbb D$. On the other hand, since $\chi_{((\alpha_i),(\beta_i))}\in (\mathbb D^\infty)^2$, $\alpha_1 \in \mathbb D$. This is a contradiction. Thus we have shown that $(\mathbb D^\infty)^2$ is closed in $\Delta(\ell^1(\mathbb N^2,\omega))$.

As $(\mathbb D^\infty)^2$ is compact and Hausdorff in both the product topology and the relative Gel'fand topology and the product topology is weaker than the relative Gel'fand topology, these topologies are equal.
\end{proof}

\begin{cor}\label{cor:6}
Let $\omega$ be a weight on $\mathbb N^2$. If $X=(\overline{B(0,r)}\times \mathbb D \times \cdots)\times \mathbb D^\infty$ is contained in $\Delta(\ell^1(\mathbb N^2,\omega))$, then the product topology and relative Gel'fand topology on $X$ are equal.
\end{cor}

Before giving a proof of Theorem \ref{thm:1} we shall require the following two lemmas. Proofs of these lemmas may be known. For completeness we shall include them.

\begin{lem}\label{lem:4}
Let $X$ be a connected topological space, and let $Y$ be a closed connected proper subset of $X$. Let $f:X \to \mathbb C$ be continuous and nowhere zero on $Y$. Then there an open subset $U$ of $X$ which properly contains $Y$ and $f(u)\neq 0$ for any $u \in U$.
\end{lem}
\begin{proof}
Take any $y \in Y$. Then there is a neighbourhood $U_y$ of $y$ such that $f(w)\neq 0$ for any $w \in U_y$. Take $U=\bigcup_{y \in Y}U_y$. Then $U$ is open in $X$ and $f(u)\neq 0$ for any $u \in U$. As $X$ is connected, $Y\subsetneq U$.
\end{proof}

\begin{lem}\label{lem:5}
Let $R>1$, and let $X=(\overline{B(0,R)}\times \mathbb D\times \mathbb D\times \cdots)\times \mathbb D^\infty$. Let $U$ be an open subset of $X$ containing $(\mathbb D^\infty)^2$. Then there is $1<r\leq R$ such that $\mathbb (D^\infty)^2\subset (\overline{B(0,r)}\times \mathbb D\times \mathbb D\times \cdots)\times \mathbb D^\infty \subset U$.
\end{lem}
\begin{proof}
Suppose that there is no $\eta>0$ such that $A_\eta=(\overline{B(0,1+\eta)}\times\mathbb D\times \cdots)\times \mathbb D^\infty$ is a subset of $U$. Then for every $n \in \mathbb N$, there is an element $Z_n=((z_1^n,z_2^n,\cdots),(w_1^n,w_2^n,\ldots))$ which is in $A_{\frac{1}{n}}$ but not in $U$ and hence not in $(\mathbb D^\infty)^2$. Since $(Z_n)$ is a sequence in a compact metric space $X$, it has a convergent subsequence. So, we may assume that $(Z_n)$ is convergent. Let $Z=((z_1,z_2,\ldots),(w_1,w_2,\ldots))$ be the limit of $(Z_n)$. Since $Z_n \to Z$, $z_i^n \to z_i$ and $w_i^n \to w_i$ as $n \to \infty$ for all $i \in \mathbb N$. Notice that $z_i^n \in \mathbb D$ for all $n \in \mathbb N\setminus\{1\}$ and $w_i^n \in \mathbb D$ for all $n \in\mathbb N$, so, $z_i \in \mathbb D$ for $i \in \mathbb N\setminus\{1\}$ and $w_i \in \mathbb D$ for all $i \in \mathbb N$. As $(Z_n)$ is a sequence in $X\setminus U$ and $U$ is open in $X$, $Z \in X\setminus U$. As $Z \notin U$ and hence not in $(\mathbb D^\infty)^2$, $z_1\notin \mathbb D$. Since $|z_1^n|\leq 1+\frac{1}{n}$ for all $n$, $|z_1|\leq 1$ and hence $z_1\in \mathbb D$. This is a contradiction. Thus there is $\eta>0$ such that $A_{\eta} \subset U$. Clearly, $A_{\eta}$ contains $(\mathbb D^\infty)^2$.
\end{proof}
The following result will be useful in proving Theorem 1. It is an analogue of Theorem 1 of \cite{HW}.
\begin{prop}\label{prop:2}
Let $\omega$ be an admissible weight on $\mathbb N^2$, and let $\widetilde a \in \mathbf A_\omega$. Then the following statements are equivalent.
\begin{enumerate}
\item[(i)] $(a(m,n))$ has inverse in $\ell^1(\mathbb N^2,\omega)$.
\item[(ii)] $\widetilde a$ is bounded away from zero.
\item[(iii)] $\widetilde a$ has inverse in $\mathbf A_\omega$.
\end{enumerate}
\end{prop}
\begin{proof}
Since $\omega$ is an admissible weight, $\Delta(\ell^1(\mathbb N^2,\omega))=\Delta(\mathbf A)=\Delta(\mathbf A_\omega)=(\mathbb D^\infty)^2$, by Corollary \ref{cor:3}. Clearly, the statements (i) and (iii) are equivalent.\\

\noindent
(ii) $\Rightarrow$ (iii). Since $\mathbf A_\omega \subset \mathbf A$, by Theorem \ref{thmA}, $\widetilde a$ has inverse in $\mathbf A$. But then $\chi(\widetilde a)\neq 0$ for all $\chi\in \Delta(\mathbf A)=\Delta(\mathbf A_\omega)$. This implies that $\widetilde a^{-1} \in \mathbf A_\omega$.\\

\noindent
(iii) $\Rightarrow$ (ii). Since $\widetilde a$ has inverse in $\mathbf A_\omega$, $\chi(\widetilde a)\neq 0$ for any $\chi \in \Delta(\mathbf A_\omega)$. As $\Delta(\mathbf A_\omega)$ is compact, there is $k>0$ such that $|\chi(\widetilde a)|\geq k$ for all $\chi \in \Delta(\mathbf A_\omega)$. For $(s_1,s_2) \in \mathbb H^2$, let $\chi_{s_1,s_2}(\widetilde b)=\sum_{(m,n)\in \mathbb N^2} b(m,n)m^{-s_1}n^{-s_2}$ for all $\widetilde b \in \mathbf A_\omega$. Then $\chi_{s_1,s_2} \in \Delta(\mathbf A_\omega)$. It follows that $|\chi_{s_1,s_2}(\widetilde a)|=|\widetilde a(s_1,s_2)|\geq k$ for all $(s_1,s_2) \in \mathbb H^2$. Therefore $\widetilde a$ is bounded away from zero.
\end{proof}

A weight $\omega:\mathbb N^2 \to [1,\infty)$ is a \emph{non-quasianalytic weight} or a \emph{Beurling-Domar weight} \cite{D} if $\sum_{n\in \mathbb N} \frac{\log(\omega(l^n,k^n))}{1+n^2}<\infty$ for all $(l,k) \in \mathbb N^2$. Note that if $\omega$ is a Beurling-Domar weight, then it is an admissible weight.

The following gives an analogue of Domar's Theorem \cite{D} in the case of Dirichlet series.
\begin{cor}\label{cor:4}
Let $\omega$ be a Beurling-Domar weight. Let $\widetilde a:\mathbb H^2 \to \mathbb C$ have $\omega$- absolutely convergent Dirichlet series. Then $\widetilde a$ has inverse in $\mathbf A_\omega$ if and only if $\widetilde a$ is bounded away from zero.
\end{cor}

\noindent
\textbf{Proof of Theorem \ref{thm:1}:}
First assume that $\omega$ is admissible. Take $\nu=\omega$ in this case. Then Proposition \ref{prop:2} implies that $\frac{1}{\widetilde a}$ has $\omega$- absolutely convergent Dirichlet series.

Now, assume that $\omega$ is not admissible. Then, by Lemma \ref{lem:2}, $\rho_i>1$ or $\mu_i>1$ for some $i$. We may assume that $\rho_1>1$. Since $\omega\geq 1$, $\mathbf A_\omega\subset \mathbf A$. As $|\widetilde a(s_1,s_2)|\geq k>0$ for all $(s_1,s_2) \in \mathbb H^2$, by Theorem \ref{thmA}, $\frac{1}{\widetilde a}$ has the representation $\frac{1}{\widetilde a}(s_1,s_2)=\sum_{(m,n)\in \mathbb N^2} b(m,n)m^{-s_1}n^{-s_2}$ for all $(s_1,s_2) \in \mathbb H^2$ and $\sum_{(m,n)\in \mathbb N^2} |b(m,n)|<\infty$. Since $\widetilde a$ is invertible in $\mathbf A$, the Banach algebra of all absolutely convergent Dirichlet series, $\widehat{\widetilde a}(\chi)\neq 0$ for every bounded semicharacter $\chi$ on $\mathbb N^2$. We know that the Gel'fand space $\Delta(\mathbf A)$ of $\mathbf A$ is homeomorphic to $(\mathbb D^\infty)^2$. As $\omega$ is an almost monotone weight, by Corollary \ref{cor:2}, the Gel'fand space of $\mathbf A_\omega$ contains $X$, where $X=(\overline{B(0,\rho_1)}\times \mathbb D\cdots)\times \mathbb D^\infty$ and, by Corollary \ref{cor:6}, the product topology on $X$ is the relative Gel'fand topology. Since $X$ contains $(\mathbb D^\infty)^2$, $\widehat{\widetilde a}$ is continuous on $X$ and $\widehat{\widetilde a}$ is nowhere zero on $(\mathbb D^\infty)^2$, by Lemma \ref{lem:4}, there is open set $U$ of $X$ containing $(\mathbb D^\infty)^2$ such that $\widehat{\widetilde a}$ is nowhere zero on $U$. By Lemma \ref{lem:5}, we get $r>1$ such that $(\overline{B(0,r)}\times \mathbb D\times \mathbb D\times \cdots)\times \mathbb D^\infty$ contains $(\mathbb D^\infty)^2$ and is contained in $U$. Define $\nu$ on $\mathbb N^2$ by $\nu(1,1)=1$, $\nu(p_1,1)=r$, $\nu(p_j,1)=1$ if $j\neq 1$ and $\nu(1,p_i)=1$ for all $i \in \mathbb N$ and extend $\nu$ on $\mathbb N^2$ so that it is multiplicative. Then $\nu$ is an almost monotone weight, $1\leq \nu \leq \omega$ and the Gel'fand space of $\mathbf A_\nu$ is $(\overline{B(0,r)}\times \mathbb D\times\mathbb D\times \cdots)\times \mathbb D^\infty$, by Lemma \ref{lem:3}. The construction $\nu$ shows that $\nu$ is constant if and only if $\omega$ is constant and that $\nu$ is admissible if and only if $\omega$ is so. Since $\widehat{\widetilde a}$ is nowhere zero on the Gel'fand space of $\mathbf A_\nu$, $\frac{1}{\widetilde a}$ is in $\mathbf A_\nu$ and hence $\frac{1}{\widetilde a}$ has $\nu$- absolutely convergent Dirichlet series.\\

\noindent
(II) The proof is inspired by \cite{BhDe}. Let $K$ be the closure of the range of $\widetilde a$. Let $(s_1,s_2) \in \mathbb H^2$. Then
\begin{eqnarray*}
|\widetilde a(s_1,s_2)|&=&\left|\sum_{(m,n)\in \mathbb N^2} a(m,n)m^{-s_1}n^{-s_2}\right|\leq \sum_{(m,n)\in \mathbb N^2} |a(m,n)|\\
&\leq & \sum_{(m,n)\in \mathbb N^2} |a(m,n)|\omega(m,n)=\|\widetilde a\|<\infty.
\end{eqnarray*} Therefore $K$ is a compact subset of $\mathbb C$. Let $\varphi$ be holomorphic on a neighbourhood $U$ of $K$. Let $\Gamma$ be a closed rectifiable Jordan contour in the open set $U$ such that $K$ is interior to $\Gamma$. Let $\lambda \in \Gamma$. Then $\lambda \notin K$. Since $K$ is compact,  $\lambda -\widetilde a$ is bounded away from zero.  Also, $\lambda -\widetilde a\in \mathbf A_\omega$. By part (I), we get a weight $\nu_\lambda$, which is constant if and only if $\omega$ is constant, such that $1\leq \nu_\lambda \leq \omega$ and $(\lambda -\widetilde a)^{-1}\in \mathbf A_{\nu_\lambda}$. Take $P_\lambda=(\lambda-\widetilde a)^{-1}$. Then $\|(\lambda-\widetilde a)^{-1}\|_{\nu_\lambda}>0$, where $\|\cdot\|_{\nu_\lambda}$ is the norm on $\ell^1(\mathbb N^2,\nu_{\lambda})$. Let $V(\lambda)=\{\mu\in \mathbb C:|\mu-\lambda|<\|P_\lambda\|_{\nu_\lambda}^{-1}\}$. Then $V(\lambda)$ is a neighbourhood of $\lambda$ and if $\mu\in V(\lambda)$, then $\mu-\widetilde a$ is invertible in $\mathbf A_{\nu_\lambda}$. For if $\mu \in V(\lambda)$, then $\mu -\widetilde a=(\lambda-\widetilde a)(1+(\mu-\lambda)P_\lambda)$ and both the elements $\lambda-\widetilde a$ and $1+(\mu-\lambda)P_\lambda$ are invertible in $\mathbf A_{\nu_\lambda}$. Since $\Gamma$ is compact, we get $\lambda_1,\ldots,\lambda_n$ in $\Gamma$ and weights $\nu_{\lambda_1},\ldots,\nu_{\lambda_n}$ such that $\Gamma\subset \bigcup_{i=1}^n V(\lambda_i)$ and if $\mu$ is in $\Gamma$, then $\mu-\widetilde a$ has inverse in $\mathbf A_{\nu_{\lambda_i}}$ for some $i$. If $\omega$ is an admissible weight, then take each of $\nu_{\lambda_i}=\omega$ and take $\xi=\omega$. In the case, when $\omega$ is not admissible, as above we may assume that $\rho_1>1$. Therefore by part (I) we get $r_i>1$, $i=1,2,\ldots,n$, such that $\nu_{\lambda_i}(p_1,1)=r_i$, $\nu_{\lambda_i}(p_j,1)=1$ if $j\neq 1$ and $\nu_{\lambda_i}(1,p_j)=1$ for all $j \in \mathbb N$. Let $r=\min\{r_1,r_2,\ldots,r_n\}$. Define $\xi$ on $\mathbb N^2$ as $\xi(p_1,1)=r$, $\xi(p_j,1)=1$ if $j\neq 1$ and $\xi(1,p_j)=1$ for all $j \in \mathbb N$ and extend it on $\mathbb N^2$ so that it is multiplicative, i.e., $\xi=\min\{\nu_{\lambda_1},\nu_{\lambda_2},\ldots,\nu_{\lambda_n}\}$. Then $\xi$ is an almost monotone weight, $\xi$ is nonconstant and $1\leq \xi\leq \omega$. Observe that $\mathbf A_\omega\subset \mathbf A_{\nu_{\lambda_i}}\subset \mathbf A_\xi$ for all $i$. Also, notice that if $\lambda \in \Gamma$, then the inverse of $\lambda -\widetilde a$ is in $\mathbf A_\xi$. Define a map $\psi:\Gamma \to \mathbf A_\xi$ by $\psi(\lambda)=\varphi(\lambda)P_\lambda$. Then $\psi$ is continuous. It  follows from functional calculus that $\varphi(\widetilde a)=\frac{1}{2\pi i}\int_\Gamma \varphi(\lambda)P_\lambda d\lambda \in \mathbf A_\xi$. Therefore $\varphi\circ \widetilde a$ has $\xi$- absolutely convergent Dirichlet series. $\hfill \Box$\\

\noindent
\textbf{Remarks.}

\noindent
\textbf{(1)} The following gives $p$-th power analogue of Theorem \ref{thm:1}.
\begin{thm}\label{thm:2}
Let $0<p\leq 1$. Let $\omega$ be an almost monotone weight on $\mathbb N^2$, and let $\widetilde a:\mathbb H^2 \to \mathbb C$ have $p$-th power $\omega$- absolutely convergent Dirichlet series. Then the following statements hold.\\
\noindent
(I) If $\widetilde a$ is bounded away from zero, then there is an almost monotone weight $\nu$ on $\mathbb N^2$ such that
\begin{enumerate}
\item $\nu$ is constant if and only if $\omega$ is constant;
\item $\nu$ is admissible if and only if $\omega$ is admissible;
\item $1\leq \nu \leq \omega$;
\item $\frac{1}{\widetilde a}$ has $p$-th power $\nu$- absolutely convergent Dirichlet series.
\end{enumerate}
(II) If $\varphi$ is holomorphic on the closure of the range of $\widetilde a$, then there is an almost monotone weight $\xi$ on $\mathbb N^2$ such that
\begin{enumerate}
\item $\xi$ is constant if and only if $\omega$ is constant;
\item $\xi$ is admissible if and only if $\omega$ is admissible;
\item $1\leq \xi \leq \omega$;
\item $\varphi \circ \widetilde a$ has $p$-th power $\xi$- absolutely convergent Dirichlet series.
\end{enumerate}
\end{thm}
To prove Theorem \ref{thm:2} we shall require the following definitions and results (Lemma \ref{lem:6}, Proposition \ref{prop:3} and Corollary \ref{cor:5}).

Let $\mathcal B$ be a complex algebra, and let $0<p\leq 1$. A map $\|\cdot\|:\mathcal B \to \mathbb R$ is a \emph{$p$-norm} on $\mathcal B$ if $\|x\|\geq 0$, $\|x\|=0$ if and only if $x=0$, $\|x+y\|\leq \|x\|+\|y\|$, $\|\alpha x\|=|\alpha|^p\|x\|$ and $\|xy\|\leq \|x\|\|y\|$ for all $x,y \in \mathcal B$ and $\alpha \in \mathbb C$. If $\mathcal B$ is complete in the $p$-norm, then $(\mathcal B,\|\cdot\|)$ is a \emph{$p$-Banach algebra} \cite{Z}.

Let $0<p\leq 1$, and let $\mathbf A_{p\omega}$ be the collection of all $\widetilde a:\mathbb H^2\to \mathbb C$ such that $$\widetilde a(s_1,s_2)=\sum_{(m,n)\in \mathbb N^2} a(m,n)m^{-s_1}n^{-s_2}\qquad((s_1,s_2)\in \mathbb H^2),$$ and  $$\|\widetilde a\|=\sum_{(m,n)\in \mathbb N^2} |a(m,n)|^p\omega(m,n)<\infty.$$
Then $\mathbf A_{p\omega}$ is a unital commutative $p$- Banach algebra and it is isometrically isomorphic to the $p$- Banach algebra $$\ell^p(\mathbb N^2,\omega)=\{a:\mathbb N^2\to \mathbb C:\|a\|=\sum_{(m,n)\in \mathbb N^2}|a(m,n)|^p\omega(m,n)<\infty\}$$ via the map $\widetilde a\mapsto a=(a(m,n))$. When $\omega=1$, we shall write $\mathbf A_p$ for $\mathbf A_{p\omega}$.

We shall use the Gel'fand theory of commutative $p$- Banach algebras developed by \.Zelazko in \cite{Z}.

\begin{lem}\label{lem:6}
Let $0<p\leq 1$, and let $\omega$ be a weight on $\mathbb N^2$. Then $\varphi$ is a complex homomorphism on $\mathbf A_{p\omega}$ if and only if $\varphi=\varphi_{\chi}$ for some $\chi \in \widehat{\mathbb N^2_\omega}$, where $$\varphi_\chi(\widetilde b)=\sum_{(m,n)\in \mathbb N^2} b(m,n)\chi(m,n)\qquad(\widetilde b \in \mathbf A_{p\omega}).$$
\end{lem}
\begin{proof}
Since $\mathbf A_{p\omega}$ and $\ell^p(\mathbb N^2,\omega)$ are isometrically isomorphic, it is sufficient to characterize complex homomorphism on $\ell^p(\mathbb N^2,\omega)$. Let $\varphi$ be a complex homomorphism on $\ell^p(\mathbb N^2,\omega)$. Since complex homomorphisms on $\ell^p(\mathbb N^2,\omega)$ are continuous, $\varphi(a)=\sum_{(m,n)\in \mathbb N^2} a(m,n)\chi(\delta_{(m,n)})$ for all $a=\sum_{(m,n)\in \mathbb N^2}a(m,n)\delta_{(m,n)}$ in $\ell^p(\mathbb N^2,\omega)$. Define $\chi:\mathbb N^2 \to \mathbb C$ by $\chi(m,n)=\varphi(\delta_{(m,n)})\;((m,n)\in \mathbb N^2)$. Then $\chi$ is a semicharacter. Also, $|\chi(m,n)|=|\varphi(\delta_{(m,n)})|\leq \|\delta_{(m,n)}\|=\omega(m,n)$ for all $(m,n)\in \mathbb N^2$. Therefore $\chi\in \widehat {\mathbb N^2}_\omega$ and $\varphi=\varphi_\chi$.

Conversely, assume that $\varphi=\varphi_\chi$ for some $\chi\in \widehat{\mathbb N^2_\omega}$. Let $a=\sum_{(m,n)\in \mathbb N^2} a(m,n)\delta_{(m,n)}$ be in $\ell^p(\mathbb N^2,\omega)$. Since $$\sum_{(m,n)\in \mathbb N^2} |a(m,n)|^p\leq \sum_{(m,n)\in \mathbb N^2} |a(m,n)|^p\omega(m,n)<\infty,$$ $a(m,n)\to 0$ as $m,n\to \infty$. Therefore there is a finite subset $F$ of $\mathbb N^2$ such that $|a(m,n)|\leq 1$ for all $(m,n)\in \mathbb N^2\setminus F$. This implies that $|a(m,n)|\leq |a(m,n)|^p$ for all $(m,n)\in \mathbb N^2\setminus F$. As a consequence $$\sum_{(m,n)\in \mathbb N^2\setminus F}|a(m,n)||\chi(m,n)|\leq \sum_{(m,n)\in \mathbb N^2\setminus F}|a(m,n)|^p\omega(m,n)$$ and hence the series $\sum_{(m,n)\in \mathbb N^2} a(m,n)\chi(m,n)$ converges absolutely. Therefore $\varphi_\chi$ is well defined and in fact it is a complex homomorphism on $\ell^p(\mathbb N^2,\omega)$.
\end{proof}

The following is an analogue of Proposition \ref{prop:2} and has a similar proof.
\begin{prop}\label{prop:3}
Let $0<p\leq 1$. Let $\omega$ be an admissible weight on $\mathbb N^2$, and let $\widetilde a \in \mathbf A_{p\omega}$. Then the following statements are equivalent.
\begin{enumerate}
\item[(i)] $(a(m,n))$ has inverse in $\ell^p(\mathbb N^2,\omega)$.
\item[(ii)] $\widetilde a$ is bounded away from zero.
\item[(iii)] $\widetilde a$ has inverse in $\mathbf A_{p\omega}$.
\end{enumerate}
\end{prop}

\begin{cor}\label{cor:5}
Let $\omega$ be a Beurling-Domar weight on $\mathbb N^2$, and let $0<p\leq 1$. Let $\widetilde a:\mathbb H^2\to \mathbb C$ have $p$-th power $\omega$- absolutely convergent Dirichlet series. Then $\widetilde a$ has inverse in $\mathbf A_{p\omega}$ if and only if $\widetilde a$ is bounded away from zero.
\end{cor}

\noindent
\textbf{Proof of Theorem \ref{thm:2}} is almost same as the proof of Theorem \ref{thm:1} with minor modifications.\\

\noindent
\textbf{(2)} Let $\widetilde a:\mathbb H \to \mathbb C$ have an $\omega$- absolutely convergent Dirichlet series and $\widetilde a$ be bounded away from zero. One would like to know whether $\frac{1}{\widetilde a}$ has $\omega$- absolutely convergent Dirichlet series. The answer is no. Define $$\widetilde a(s)=2+2^{-s} \qquad(s \in \mathbb H).$$ We define a weight $\omega$ as follows. Let $n \in \mathbb N$. Then it has unique prime factorization $n=2^kp_1^{k_1}\cdots p_l^{k_l}$, where $k, l\in \mathbb N_0$, $k_i \in \mathbb N$, $p_i$'s are distinct primes and $p_i\neq 2$ for any $i$. Define $\omega(n)=2^{k+1}$. Clearly, $\widetilde a$ has $\omega$- absolutely convergent Dirichlet series. We now show that $\widetilde a$ is bounded away from zero. Let $s=\sigma+it \in \mathbb H$. Then $\sigma \geq 0$. Now,
\begin{eqnarray*}
|\widetilde a(s)|^2 & = & |2+2^{-s}|^2=|2+2^{-\sigma}e^{-it\log 2}|^2\\
& =& (2+2^{-\sigma}\cos(t\log 2))^2+\sin^2(t\log 2)\\
& = & 4+4\cdot2^{-\sigma}\cos(t\log 2)+2^{-2\sigma}\\
&\geq & 4-4\cdot2^{-\sigma}+2^{-2\sigma} = (2-2^{-\sigma})^2.
\end{eqnarray*}
Since $\sigma\geq 0$, $2^{-\sigma}\leq 1$. Therefore $|\widetilde a(s)|\geq 1$ for all $s \in \mathbb H$, i.e., $\widetilde a$ is bounded away from zero. An easy computation shows that $$\frac{1}{\widetilde a}(s)=\sum_{n=1}^\infty \frac{(-1)^n}{2^{n+1}}(2^n)^{-s}\qquad(s\in \mathbb H).$$ Clearly, $\frac{1}{\widetilde a}$ does not have $\omega$- absolutely convergent Dirichlet series.\\

\noindent
\textbf{(3)} A proof for a function of $d$ variables can naturally be given using some modifications as below. A weight $\omega$ on $\mathbb N^d$ is almost monotone if either $\omega$ is admissible, i.e., $\lim_{n\to \infty}\omega(m_1^n,\ldots,m_d^n)^{\frac{1}{n}}=1\;((m_1,\ldots,m_d)\in \mathbb N^d)$, or there is $K>0$ such that $\omega(m_1,\ldots,m_d)\leq \omega(n_1,\ldots,n_d)$ whenever $(m_1,\ldots,m_d), (n_1,\ldots,n_d)\in \mathbb N^d$ and $m_i|n_i$ for some $i$. For $j=1,\ldots,d$ and $i\in \mathbb N$, let $$\rho_i^j=\lim_{n\to \infty}\omega(1,1,\ldots,1,p_i^n,1,\ldots,1)^{\frac{1}{n}},$$ where $p_i$ is at the $j$-th place. Now, proofs for a function of $d$- variables will be simple modifications of the proofs given here.\\

\noindent
\textbf{(4)} We do not know whether there is a weight $\omega$ on $\mathbb N^2$ or $\mathbb N$ which is not an almost monotone weight.\\

\noindent
\textbf{(5)} The Gel'fand space of $\ell^1(\mathbb N^2,\omega)$, for an arbitrary weight $\omega$ on $\mathbb N^2$, is not known to us. Here, we have used a part of it since it sufficient for our purpose. It seems interesting to identify the Gel'fand space of $\ell^1(\mathbb N^2,\omega)$.



\begin{thebibliography}{33}
\bibitem{BDD} S. J. Bhatt, P. A. Dabhi and H. V. Dedania, \emph{Beurling algebra analogues of theorems of Wiener - L\'evy - \.Zelazko and \.Zelazko}, Studia Mathematica, \textbf{195}(3)(2009), 219 -- 225.
\bibitem{BhDe} S. J. Bhatt and H. V. Dedania, \emph{Beurling algebra analogues of the classical theorems of Wiener and L\'evy on absolutely convergent Fourier series}, Proc. Indian Acad. Sci. (Math. Sci.), \textbf{113}(2)(2003), 179 -- 182.
\bibitem{D} Y. Domar, \emph{Harmonic analysis based on certain commutative Banach algebras}, Acta Math., \textbf{96}(1956), 1 -- 66.
\bibitem{Ed} D. A. Edwards, \emph{On absolutely convergent Dirichlet series}, Proc. Amer. Math. Soc., \textbf{8}(1957), 1067 -- 1074.
\bibitem{GRS} I. M. Gel'fand, D. Ra\v ikov, G. E. \v Silov, \emph{Commutative normed rings}, Chelsea Publishing Company, New York, 1964.
\bibitem{GL} Helge Gl\"ockner and Lutz G. Lucht, \emph{Weighted inversion of general Dirichlet series}, Trans. Amer. Math. Soc., \textbf{366}(6)(2014), 3275 -- 3293.
\bibitem{GN} Arthur Goodman and D. J. Newman, \emph{A Wiener type theorem for Dirichlet series}, Proc. Amer. Math. Soc., \textbf{92}(4)(1984), 521 -- 527.
\bibitem{HW} Edwin Hewitt and J. H. Williamson, \emph{Note on absolutely convergent Dirichlet series}, Proc. Amer. Math. Soc., \textbf{8}(1957), 863 -- 868.
\bibitem{HZ} Edwin Hewitt and H. S. Zuckermann, \emph{The $\ell_1$-algebra of a commutative semigroup}, Trans. Amer. Math. Soc., \textbf{83}(1)(1956), 70 -- 97.
\bibitem{L} P. L\'evy, \emph{Sur la convergence absolue des s\'eries de Fourier}, Composito Mathematica, \textbf{1}(1935), 1 -- 14.
\bibitem{CR} C. E. Rickart, \emph{General Theory of Banach Algebras}, D. Van Nostrand Inc., Princeton, New Jersey, 1960.
\bibitem{W} N. Wiener, Tauberian theorems, \emph{Annals of Mathematics}, \textbf{33}(1932), 1 -- 100.
\bibitem{Z} W. \.Zelazko, \emph{Selected topics in topological algebras}, Aarhus Universitet, Lecture Notes Series No. 31, 1971.
\end{thebibliography}
\end{document}